\documentclass [11pt] {article}
\usepackage{amsfonts}
\usepackage{amssymb}
\usepackage{theorem}

\newtheorem{lemma}{Lemma}[section]

\newtheorem{theorem}{Theorem}

\newtheorem{definition}{Definition}[section]

\begin{document}

\title{The Symmetry Preserving Removal Lemma}

\author{\sc Bal\'azs Szegedy}

\maketitle

\begin{abstract} In this note we observe that in the hyper-graph removal lemma the edge removal can be done in a way that the symmetries of the original hyper-graph remain preserved. As an application we prove the following generalization of Szemer\'edi's Theorem on arithmetic progressions. If in an Abelian group $A$ there are sets $S_1,S_2\dots,S_t$ such that the number of arithmetic progressions $x_1,x_2,\dots,x_t$ with $x_i\in S_i$ is $o(|A|^2)$ then we can shrink each $S_i$ by $o(|A|)$ elements such that the new sets don't have such a diagonal arithmetic progression.
\end{abstract}

\section{Introduction}

A directed $k$-uniform hyper-graph $H$ on the vertex set $V$ is a subset of $V^k$ such that there is no repetition in the $k$ coordinates. A {\bf homomorphism} between two directed $k$-uniform hyper-graphs $F$ and $H$ with vertex sets $V(F)$ and $V(H)$ is a map $f:V(F)\mapsto V(H)$  such that $(f(a_1),f(a_2),\dots f(a_k))$ is in $H$ whenever $(a_1,a_2,\dots a_k)$ is in $F$. The {\bf automorphism group} ${\rm Aut}(H)$ is the group of bijective homomorphisms $\pi:V(H)\mapsto V(H)$.
The {\bf homomorphism density} $t(F,G)$ of $F$ in $G$ is the probability that a random map $f:V(G)\mapsto V(H)$ is a homomorphism.

The so-called hyper-graph removal lemma (\cite{NRS},\cite{RS},\cite{Gow},\cite{Ish},\cite{Tao})(in the directed setting) says the following

\begin{theorem}[Removal Lemma]\label{rem} For every $k\in\mathbb{N}$, $\epsilon>0$ and $k$-uniform directed hyper-graph $F$ there is a constant $\delta=\delta(k,\epsilon,F)>0$ such that for every $k$-uniform directed hyper-graph $G$ with $t(F,G)\leq\delta$ there is a subset $S\subseteq G$ with $S\leq\epsilon |V(G)|^k$ such that $t(F,G\setminus S)=0$.
\end{theorem}

Using this deep result we observe that the edge removal can be done in a way that the symmetries of $G$ remain preserved.

\begin{theorem}[Symmetry Preserving Removal Lemma]\label{symprem} For every $k\in\mathbb{N}$, $\epsilon>0$ and $k$-uniform directed hyper-graph $F$ there is a constant $\delta_2=\delta_2(k,\epsilon,F)>0$ such that for every $k$-uniform directed hyper-graph $G$ with $t(F,G)\leq\delta_2$ there is a subset $S\subseteq G$ with $S\leq\epsilon |V(G)|^k$ such that $t(F,G\setminus S)=0$ and furthermore ${\rm Aut}(G)\subseteq {\rm Aut}(G'\setminus S)$.
\end{theorem}

\begin{proof} Let $V=V(G)$. Using the original removal lemma it remains to show that if $S\subseteq V^k$ satisfies $t(F,G\setminus S)=0$ than there is $S'\subseteq V^k$ which is ${\rm Aut}(G)$ invariant, $t(F,G\setminus S')=0$ and $S'\leq |F||S|$. Such an $S'$ is the union of those orbits $O$ of ${\rm Aut}(G)$ on $V^k$ for which $|O|/|O\cap S|\leq |F|$. Assume that $f:V(F)\mapsto V$ is a homomorphism from $F$ to $G\setminus S'$. Then for every fixed $e\in F$ and for random element $\pi\in{\rm Aut}(G)$ the probability that $\pi(f(e))\in G\setminus S$ is less that $1/|F|$ and so there is some $\pi\in\rm{Aut}(G)$ with $\pi(f(F))\subseteq G\setminus S$ which is contradiction.
\end{proof}

The argument given for the symmetry preserving removal is very general. It applies for various modified versions of the removal lemma. An important such version is the $t$-partite removal lemma where $t$ is a fixed natural number. The vertex set of a $t$-partite $k$-uniform hypergraph is a $t$ tuple $\{V_i\}_{i=1}^t$ of finite sets. An edge of a $t$-partite hypergraph is an element from $\prod_{i=1}^kV_{a_i}$ where $a_1,a_2,\dots,a_k$ are $k$ distinct numbers between $1$ and $t$.
Let $G_1,G_2$ be two $t$ partite $k$-uniform hyper-graphs with vertex sets $\{V_i\}_{i=1}^t$ and $\{W_i\}_{i=1}^t$. A homomorphism from $G_1$ to $G_2$ is a $t$ tuple of maps $\{\phi_i:V_i\to W_i\}_{i=1}^t$ such that $(\phi_{a_1}(r_1),\phi_{a_2}(r_2),\dots,\phi_{a_k}(r_k))\in\prod_{i=1}^k W_{a_i}$ is an edge in $G_2$ whenever $(r_1,r_2,\dots,r_k)\in\prod_{i=1}^k V_{a_i}$ is an edge in $G_1$. An automorphism is a bijective homomorphism from $G_1$ to $G_1$ and the homomorphism density $t(G_1,G_2)$ is the probability that a random $t$ tuple of maps $\{\phi_i:V_i\to W_i\}_{i=1}^t$ is a homomorphism.

We give an example for an application of the symmetry preserving removal lemma and then we generalize it in the next chapter.
\bigskip

\noindent {\bf Example 1.:} Let $S$ be a subset of a group $T$. The Cayley graph ${\rm Cy}(T,S)\subseteq G\times G$ is the collection of pairs $(a,b)$ with $ab^{-1}\in S$. The automorphism group of ${\rm Cy}(T,S)$ contains $T$ with the action $(a,b)\mapsto (ag,bg)$. Clearly any subset of $T\times T$ invariant under this action of $G$ is a Cayley graph itself. This means that the $T$-orbit of edges in ${\rm Cy(T,S)}$ correspond to elements of $S$.
 We apply the symmetry preserving removal lemma for $F=\{(1,2),(1,3),(2,3)\}$ with $V(F)=\{1,2,3\}$ and for $G={\rm Cy}(T,S)$. A homomorphism from $F$ to $G$ is a map $f:\{1,2,3\}\mapsto T$ such that $a=f(1)f(2)^{-1}$ , $b=f(2)f(3)^{-1}$ and $c=f(1)f(3)^{-1}$ are all in $S$. Consequently, the number of such homomorphisms is the number is $|T||\{(a,b,c)|ab=c~,~a,b,c\in S\}$.
 The symmetry preserving removal lemma yields that if $ab=c$ has $o(|T|^2)$ solutions in $S$ then one can remove a $o(|T|)$ elements from $S$ such that in the new set there is no solution of $ab=c$. This was first proved by Ben Green \cite{Green} for Abelian groups and generalized for groups by Kral, Serra and Vena \cite{ksl}.

\section{Cayley Hypergraphs}

In this chapter we describe a potential way of generalizing Cayley graphs to the hypergraphs setting and then discuss the symmetry preserving removal lemma on such graphs.

\begin{definition} Let $G_1,G_2,\dots,G_t$ be $t$ finite groups and let $H$ be a subgroup of $\prod_{i=1}^t G_i$. The group $H$ is acting on each $G_i$ by $(h_1,h_2,\dots,h_t)g=h_ig$ where $(h_1,h_2,\dots,h_t)\in H$ and $g\in G_i$. A $t$ partite $k$-uniform hypergraph $T$ on the vertex set $\{G_i\}_{i=1}^t$ is called a Cayley hypergraph if its automorphism group contains $H$ with the previous action.
\end{definition}

This definition is very general so we will start to analyze a special setting.
Assume that all the groups $G_1,G_2,\dots,G_t$ are isomorphic to an Abelian group $A$. Furthermore, to get something interesting we want to assume that $H$ is not too big and not too small. Let $C=\{C_1,C_2,\dots,C_r\}$ be a collection of $k$-element subsets of $\{1,2,\dots,t\}$. Each set $C_i$ defines a projection $p_i:H\mapsto A^k$ to the coordinates in $C_i$. Assume that the factor group $A^{C_i}/p_i(H)\cong A$ and let $\psi_i:A^{C_i}\to A$ be a homomorphism with kernel $p_i(H)$. Now we pick a subsets $S_i\subseteq A$ for $1\leq i\leq r$ and we define the graph $H_{k,t}(A,\{S_i\},C)$ as
$$\bigcup_{i=1}^r \psi_i^{-1}(S_i)$$
where $\psi_i^{-1}(S_i)$ is the union of cosets in $A^{C_i}$ of $p_i(H)$ representing an element in $S_i$. Note that the way we produced $H_{k,t}(A,\{S_i\},C)$ guarantees that its automorphism group contains $H$ as a subgroup.

The symmetry preserving removal lemma for $t$-partite hypergraphs directly implies the following lemma:

\begin{lemma}[Cayley Hypergraph Removal Lemma]\label{CHRL} For every $k,t$ natural numbers and $\epsilon>0$ there exists a constant $\delta>0$ such that if $$t(F,H_{k,t}(A,\{S_i\},C))\leq\delta$$ for some $t$ partite $k$-uniform hypergraph $F$ then there are subsets $S_i'$ in $A$ of size at most $\epsilon|A|$ such that $t(F,H_{k,t}(A,\{S_i\setminus S_i'\},C))=0$.
\end{lemma}

\bigskip


\noindent{\bf Example 2.:} This example uses an idea by Solymosi \cite{S} who showed that the Hypergraph Removal Lemma implies Szemere'di's theorem on arithmetic progressions (even in a multi dimensional setting). Let $t$ be a natural number, $k=t-1$ and $A$ be an Abelian group. We define $H$ to be the subgroup in $A^t$ of the elements $(a_1,a_2,\dots,a_t)$ with $\sum_{i=1}^ta_i=0$ and $\sum_{i=1}^t(i-1)a_i=0$.
Now set $$C=\{\{1,2,\dots,i-1,i+i,\dots,t\}\}_{i=1}^t$$ and
$$\psi_i(a_1,a_2,\dots,a_{i-1},a_{i+1},\dots,a_t)=\sum_{j=1}^t (j-i)a_i.$$
The functions $\psi_i$ are computed in a way that ${\rm ker}(\psi_i^{-1})$ is the projection of $H$ to the coordinates in $C_i$.

Let $F$ be the complete $t$ partite $t-1$ uniform hypergraph on four point. Lemma $\ref{CHRL}$ applied to $F$ and the above hypergraph $H_{t-1,t}(A,\{S_i\},C)$ implies that if the system
$$x_i=\sum_{j=1}^t (j-i)a_i\in S_i$$ has $o(|A|^t)$ solutions then we can delete $o(|A|)$ elements from each $S_i$ such that the previous system has no solution. It is clear that $x_1,x_2,\dots,x_t$ are forming a $t$ term arithmetic progression and in fact any such progression with $x_i\in S_i$ gives rise to $|A|^{t-2}$ solution of the previous system. Using this we obtain the following:

\medskip

\begin{theorem}[Diagonal Szemer\'edi Theorem] For every $\epsilon>0$ there exists a $\delta>0$ such that if $A$ is an Abelian group, $S_1,S_2,\dots,S_t$ are subsets in $A$ and there are at most $\delta|A|^2$ $t$-tuples $x_1,x_2,\dots,x_t$ with $x_i\in S_i$ such that they are forming a $t$ term arithmetic progression then we can shrink each $S_i$ by at most $\epsilon|A|$ elements such that the new sets don't have such a configuration.
\end{theorem}

This theorem implies Szemer\'edi's theorem \cite{Szem} if we apply it for $S=S_i,~1\leq i\leq t$ since $S$ contains the trivial progressions $a,a,a\dots,a$ which are only removable if we delete the whole set $S$.

\end{document}